\documentclass[a4paper,12pt]{article}
\usepackage{amssymb}
\usepackage[latin1]{inputenc}     

\newcommand{\beq}{\begin{equation} }
\newcommand{\eqq}{\end{equation} }
\newcommand{\cuad}{{\sqcap\kern-.68em\sqcup}}

\newcommand{\equ}[1]{(\ref{#1})}

\newtheorem{remark}{Remark}[section]
\newcommand{\bremark}{\begin{remark} \em}
\newcommand{\eremark}{\end{remark} }

\def\beeq{\begin{equation}}
\def\eeq{\end{equation}}
\newcommand{\begeqaet}{\begin{eqnarray*}}
\newcommand{\eneqaet}{\end{eqnarray*}}

\let\Section=\section
\def\section{\setcounter{equation}{0}\Section}
\newtheorem{Lem}{Lemma}[section]
\newtheorem{Thm}{Theorem}[section]

\newtheorem{Remark}{Remark}[section]

\begin{document}
\begin{center}{\bf\Large Ground state solution for  a class of differential equations with left and right fractional derivatives}\medskip

\bigskip

\bigskip

{C\'esar Torres}

 Departamento de Ingenier\'{\i}a  Matem\'atica and
Centro de Modelamiento Matem\'atico
 UMR2071 CNRS-UChile,
 Universidad de Chile\\
 Santiago, Chile.\\
 {\sl  (ctorres@dim.uchile.cl)}

\end{center}

\medskip

\medskip
\medskip
\medskip
\medskip

\begin{abstract}
In this work we study the existence of solution for a class of fractional differential equation given by
\begin{eqnarray}\label{eq00}
_{t}D_{\infty}^{\alpha}{_{-\infty}}D_{t}^{\alpha}u(t) + u(t)  = & f(t,u(t))\\
u\in H^{\alpha}(\mathbb{R}).\nonumber
\end{eqnarray}
where $\alpha \in (1/2, 1)$, $t\in \mathbb{R}$, $u\in \mathbb{R}$, $f\in C(\mathbb{R}, \mathbb{R})$. Using mountain pass theorem and comparison argument we prove that (\ref{eq00}) at least has one nontrivial solution.
\end{abstract}
\date{}

\setcounter{equation}{0}
\section{ Introduction}

The aim of this article is to study the fractional differential equation with left and right fractional derivative
\begin{eqnarray}\label{Eq00}
_{t}D_{\infty}^{\alpha}{_{-\infty}}D_{t}^{\alpha}u(t) + u(t)  = & f(t,u(t))\\
u\in H^{\alpha}(\mathbb{R}).\nonumber
\end{eqnarray}
where $\alpha \in (1/2, 1)$, $t\in \mathbb{R}$, $u\in \mathbb{R}$, $f\in C(\mathbb{R}, \mathbb{R})$.

The study of fractional calculus (differentiation and integration of arbitrary order) has emerged as an important and popular field of research. It is mainly due to the extensive application of fractional differential equations in many engineering and scientific disciplines such as physics, chemistry, biology, economics, control theory, signal and image processing, biophysics, blood flow phenomena, aerodynamics, fitting of experimental data, etc., \cite{OAJTMJS}, \cite{RH}, \cite{AKHSJT}, \cite{KMBR}, \cite{IP}, \cite{JSOAJTM}, \cite{GZ}. An important characteristic of fractional-order differential operator that distinguishes it from the integer-order differential operator is its nonlocal behavior, that is, the future state of a dynamical system or process involving fractional derivative depends on its current state as well its past states. In other words, differential equations of arbitrary order describe memory and hereditary properties of various materials and processes. This is one of the features that has contributed to the popularity of the subject and has motivated the researchers to focus on fractional order models, which are more realistic and practical than the classical integer-order models.

Very recently, also equations including both left and right fractional derivatives were investigated \cite{TABS-1}, \cite{TABS-2}, \cite{MK-1}, \cite{MK-2}, \cite{MTDB}. Such differential equations mixing both types of derivatives found interesting applications in fractional variational principles, fractional control theory, as well as fractional Lagrangian and Hamiltonian dynamics \cite{OA}, \cite{DB}, \cite{DBSM}, \cite{DBSMKT}, \cite{ERKNRHSMDB}. 

Although investigations concerning ordinary and partial fractional differential equations yield many interesting and important results (compare that enclosed in monographies, \cite{AKHSJT}, \cite{KMBR}, \cite{IP}, \cite{SSAKOM}) for equations with operators including fractional derivatives of one type, still the fractional equations with mixed derivatives need further study: Using the fixed point theorems \cite{OA1}, \cite{DBJT}, \cite{MK-0} one can obtain analytical results, this solution is represented by series of alternately left and right fractional integrals. In \cite{MK-00}, Klimek shows an application of the Mellin transform, but this solution is represented by complicated series of special functions.

It should be noted that critical point theory and variational methods have also turned out to be very effective tools in determining the existence of solutions for integer order differential equations. The idea behind them is trying to find solutions of a given boundary value problem by looking for critical points of a suitable energy functional defined on an appropriate function space. In the last 30 years, the critical point theory has become to a wonderful tool in studying the existence of solutions to differential equations with variational structures, we refer the reader to the books due to Mawhin and Willem \cite{JMMW}, Rabinowitz \cite{PR} and the references listed therein. 

Motivated by the above classical works, in recent paper \cite{FJYZ}, for the first time, the authors showed that the critical point theory is an effective approach to tackle the existence of solutions for the following fractional boundary value problem
\begin{eqnarray}\label{Eq01}
&_{t}D_{T}^{\alpha}({_{0}D_{t}^{\alpha}}u(t)) = \nabla F(t,u(t)),\;\mbox{a.e.}\;t\in [0,T],\\
&u(0) = u(T) = 0.\nonumber
\end{eqnarray}
and obtained the existence of at least one nontrivial solution. In \cite{CT-1} using the same tool, the author studied the following fractional boundary value problem
\begin{eqnarray}\label{Eq02}
&_{t}D_{T}^{\alpha}({_{0}D_{t}^{\alpha}}u(t)) = f(t,u(t)),\;\;\;t\in [0,T],\\
&u(0) = u(T) = 0.\nonumber
\end{eqnarray} 
where $\frac{1}{2}< \alpha < 1$ and $f\in C(\mathbb{R}\times \mathbb{R})$ is super-quadratic at infinity and sub-quadratic at the origin. We note that it is not easy to use the critical point theory to study (\ref{Eq01}) and (\ref{Eq02}), since it is often very difficult to establish a suitable space and variational functional for the fractional boundary value problem. 

Inspired by theses works, in \cite{CT-2} the author considered the following fractional systems
\begin{eqnarray}\label{Eq03}
&_{t}D_{\infty}^{\alpha}({_{-\infty}D_{t}^{\alpha}}u(t)) + L(t)u(t)= \nabla W(t,u(t)),\\
&u \in H^{\alpha}(\mathbb{R}, \mathbb{R}^{n}),\nonumber
\end{eqnarray} 
where $\alpha \in (1/2,1)$, $t\in \mathbb{R}$, $u\in \mathbb{R}^{n}$, $L\in C(\mathbb{R}, \mathbb{R}^{n\times n})$ is a symmetric matrix valued function and $W:\mathbb{R}\times \mathbb{R}^{n} \to \mathbb{R}$; satisfies the following condition
\begin{itemize}
\item[$(L)$] $L(t)$ is positive definite symmetric matrix for all $t\in \mathbb{R}$ and there exists an $l\in C(\mathbb{R}, (0,\infty))$ such that $l(t) \to +\infty$ as $t \to \infty$ and
    \begin{equation}\label{Eq04}
    (L(t)x,x) \geq l(t)|x|^{2},\;\;\mbox{for all}\;t\in \mathbb{R}\;\;\mbox{and}\;\;x\in \mathbb{R}^{n}.
    \end{equation}
\item[$(W_{1})$] $W\in C^{1}(\mathbb{R} \times \mathbb{R}^{n}, \mathbb{R})$ and there is a constant $\mu >2$ such that
$$
0< \mu W(t,x) \leq (x, \nabla W(t,x)),\;\;\mbox{for all}\;t\in \mathbb{R}\;\;\mbox{and}\;x\in\mathbb{R}^{n}\setminus \{0\}.
$$
\end{itemize}
Assuming that (L) and ($W_{1}$) hold and introducing some compact embedding theorem, in \cite{CT-2} the author showed that the (PS) condition is satisfied and obtained the existence of solutions of (\ref{Eq03}) using the usual Mountain Pass Theorem. However the ($L$) condition does not seem to be natural and is restrictive. For example, if $L(t) = sI_{n}$ (where $s>0$ is a constant and $I_{n}$ is the unit matrix of order $n$), then (\ref{Eq04}) does not cover this trivial case. In \cite{CT-3}, the author showed that the condition (\ref{Eq04}) can be removed if $L$ is uniformly bounded from below, that is
\begin{itemize}
\item[($L_{1}$)] $L\in C(\mathbb{R}, \mathbb{R}^{n^{2}})$ is a symmetric and positive definite matrix for all $t\in \mathbb{R}$ and there is a constant $M>0$ such that
$$
(L(t)x,x)\geq M|x|^{2}\;\;\mbox{for}\;\;t\in \mathbb{R}\;\mbox{and}\;x\in \mathbb{R}^{n},
$$
\end{itemize}
and $W$ has the form $W(t,u) = a(t)V(u)$, where
\begin{itemize}
\item[($V_{1}$)] $a:\mathbb{R} \to \mathbb{R}^{+}$ is a continuous function such that 
$$
\lim_{|t|\to +\infty} a(t) = 0;
$$
\item[($V_{2}$)] $V\in C^{1}(\mathbb{R}^{n}, \mathbb{R})$ and there is a constant $\mu >2$ such that
$$
0< \mu V(u) \leq (\nabla V(u), u)\;\;\mbox{for all}\;\;u\in \mathbb{R}^{n} \setminus \{0\} .
$$
\end{itemize}
Under $(L_{1}) - (V_{2})$, the author obtained the existence solutions of (\ref{Eq03}) using again the Mountain Pass Theorem.

In this article we are interested in studying the nonlinear fractional differential equations with left and right fractional derivatives given by (\ref{Eq00}). Our goal is to study the existence of ground states of (\ref{Eq00}). Now we state our main assumptions. In order to find solutions of (\ref{Eq00}), we will assume the following general hypotheses.

\begin{itemize}
\item[($f_{0}$)] $f(t, \xi) \geq  0$ if $\xi \geq 0$ and $f(t, \xi) = 0$ if $\xi \leq 0$, for all $t\in \mathbb{R}$ 
\item[($f_{1}$)] There exists $\theta >2$ such that
$$
0 < \theta F(t,\xi) \leq \xi f(t, \xi),\;\;\forall (t,\xi),\;\xi \neq 0,
$$ 
where $F(t,\xi) = \int_{0}^{\xi} f(t, \sigma)d\sigma$.
\item[$(f_{2})$]  $f(t,\xi) = o(|\xi|)$ uniformly in $t$.
\item[($f_{3}$)] $\lim_{|\xi|\to \infty} \frac{f(t,\xi)}{|\xi|^{p_{0}}} = 0$ for some $p_{0} + 1 > \theta$, uniformly in $t\in \mathbb{R}$.
\item[$(f_{4})$] $\frac{f(t, \sigma \xi) \xi}{\sigma}$ is a increasing function for every $\sigma >0$, $t, \xi \in \mathbb{R}$. 
\item[$(f_{5})$] There exist continuous function $\overline{f}$ and $a$, defined in $\mathbb{R}$, such that $\overline{f}$ satisfies $(f_{0})-(f_{4})$  and
\begin{eqnarray*}
&0 \leq f(t, \xi) - \overline{f}(\xi) \leq a(t) (|\xi| + |\xi|^{p_{0}})\quad \mbox{for all}\;\;t, \xi \in \mathbb{R}, \\
& \lim_{|t| \to \infty} a(t) = 0
\end{eqnarray*}
and
$$
m(\{t\in \mathbb{R}: f(t,\xi) > \overline{f}(\xi)\})>0,
$$
where $m$ denotes the Lebesgue measure.
\end{itemize}

At this point we state our existence theorem for the autonomous equation, that is for $\overline{f}$. This theorem will server as a basis for the proof of the main existence theorem for the case where $f$ depends on $x$. 

\begin{Thm}\label{tm01}
Assume that $\frac{1}{2}< \alpha < 1$ and that $\overline{f}:\mathbb{R} \to \mathbb{R}$ is a function. Then we have the following. If $\overline{f}$ satisfies $(f_{0})-(f_{4})$, then 
$$
{_{t}}D_{\infty}^{\alpha} {_{-\infty}}D_{t}^{\alpha}u(t) + u(t) = \overline{f}(u)\;\;\mbox{in }\;\;\mathbb{R}
$$
has a nontrivial weak solution.
\end{Thm}

Now we state our main existence theorem

\begin{Thm}\label{tm02}
Assume that $\frac{1}{2}< \alpha < 1$ and $f$ satisfies $(f_{0})-(f_{5})$. Then equation (\ref{Eq00}) possesses at least one weak nontrivial solution. 
\end{Thm}

We prove the existence of weak solution of (\ref{Eq00}) applying the mountain pass theorem \cite{AAPR} to the functional $I$ defined on $H^{\alpha}(\mathbb{R})$ as 
\begin{equation}\label{Eq06}
I(u) = \frac{1}{2}\int_{\mathbb{R}} [|{_{-\infty}}D_{t}^{\alpha}u(t)|^{2} + u(t)^{2}]dt - \int_{\mathbb{R}}F(t, u(t))dt.
\end{equation}
However, the direct application of the mountain pass theorem is not sufficient, since the Palais-Smale sequence might lose compactness in the whole space $\mathbb{R}$. To overcome this difficulty, we use a comparison argument devised in \cite{PCOM} and \cite{PR1} for $\alpha = 1$ and by \cite{PFAQJT} for the fractional laplacian, based on the energy functional 
\begin{equation}\label{Eq07}
\overline{I}(u) = \frac{1}{2} \int_{\mathbb{R}} [|{_{-\infty}}D_{t}^{\alpha}u(t)|^{2} + u(t)^{2}]dt - \int_{\mathbb{R}}\overline{F}(u(t))dt.
\end{equation} 

The rest of the paper is organized as follows: in section 2, subsection 2.1, we describe the Liouville-Weyl fractional calculus; in subsection 2.2 we introduce the fractional space that we use in our work and some proposition are proven which will aid in our analysis. In section 3, we will prove Theorem \ref{tm01} and Theorem \ref{tm02}.
\section{Preliminary Results}

\subsection{Liouville-Weyl Fractional Calculus}

The Liouville-Weyl fractional integrals of order $0<\alpha < 1$ are defined as
\begin{equation}\label{LWeq01}
_{-\infty}I_{x}^{\alpha}u(x) = \frac{1}{\Gamma (\alpha)} \int_{-\infty}^{x}(x-\xi)^{\alpha - 1}u(\xi)d\xi
\end{equation}
\begin{equation}\label{LWeq02}
_{x}I_{\infty}^{\alpha}u(x) = \frac{1}{\Gamma (\alpha)} \int_{x}^{\infty}(\xi - x)^{\alpha - 1}u(\xi)d\xi
\end{equation}
The Liouville-Weyl fractional derivative of order $0<\alpha <1$ are defined as the left-inverse operators of the corresponding Liouville-Weyl fractional integrals
\begin{equation}\label{LWeq03}
_{-\infty}D_{x}^{\alpha}u(x) = \frac{d }{d x} {_{-\infty}}I_{x}^{1-\alpha}u(x)
\end{equation}
\begin{equation}\label{LWeq04}
_{x}D_{\infty}^{\alpha}u(x) = -\frac{d }{d x} {_{x}}I_{\infty}^{1-\alpha}u(x)
\end{equation}
The definitions (\ref{LWeq03}) and (\ref{LWeq04}) may be written in an alternative form:
\begin{equation}\label{LWeq05}
_{-\infty}D_{x}^{\alpha}u(x) = \frac{\alpha}{\Gamma (1-\alpha)} \int_{0}^{\infty}\frac{u(x) - u(x-\xi)}{\xi^{\alpha + 1}}d\xi
\end{equation}
\begin{equation}\label{LWeq05}
_{x}D_{\infty}^{\alpha}u(x) = \frac{\alpha}{\Gamma (1-\alpha)} \int_{0}^{\infty}\frac{u(x) - u(x+\xi)}{\xi^{\alpha + 1}}d\xi
\end{equation}

\noindent
We establish the Fourier transform properties of the fractional integral and fractional differential operators. Recall that the Fourier transform $\widehat{u}(w)$ of $u(x)$ is defined by
$$
\widehat{u}(w) = \int_{-\infty}^{\infty} e^{-ix.w}u(x)dx.
$$
Let $u(x)$ be defined on $(-\infty, \infty)$. Then the Fourier transform of the Liouville-Weyl integral and differential operator satisfies
\begin{equation}\label{LWeq06}
\widehat{ _{-\infty}I_{x}^{\alpha}u(x)}(w) = (iw)^{-\alpha}\widehat{u}(w)
\end{equation}
\begin{equation}\label{LWeq07}
\widehat{ _{x}I_{\infty}^{\alpha}u(x)}(w) = (-iw)^{-\alpha}\widehat{u}(w)
\end{equation}
\begin{equation}\label{LWeq08}
\widehat{ _{-\infty}D_{x}^{\alpha}u(x)}(w) = (iw)^{\alpha}\widehat{u}(w)
\end{equation}
\begin{equation}\label{LWeq09}
\widehat{ _{x}D_{\infty}^{\alpha}u(x)}(w) = (-iw)^{\alpha}\widehat{u}(w)
\end{equation}
\subsection{Fractional Derivative Spaces}

In this section we introduce some fractional spaces for more detail see \cite{VEJR}.

\noindent
Let $\alpha > 0$. Define the semi-norm
$$
|u|_{I_{-\infty}^{\alpha}} = \|_{-\infty}D_{x}^{\alpha}u\|_{L^{2}}
$$
and norm
\begin{equation}\label{FDEeq01}
\|u\|_{I_{-\infty}^{\alpha}} = \left( \|u\|_{L^{2}}^{2} + |u|_{I_{-\infty}^{\alpha}}^{2} \right)^{1/2},
\end{equation}
and let
$$
I_{-\infty}^{\alpha} (\mathbb{R}) = \overline{C_{0}^{\infty}(\mathbb{R})}^{\|.\|_{I_{-\infty}^{\alpha}}}.
$$
Now we define the fractional Sobolev space $H^{\alpha}(\mathbb{R})$ in terms of the fourier transform. Let $0< \alpha < 1$, let the semi-norm
\begin{equation}\label{FDEeq02}
|u|_{\alpha} = \||w|^{\alpha}\widehat{u}\|_{L^{2}}
\end{equation}
and norm
$$
\|u\|_{\alpha} = \left( \|u\|_{L^{2}}^{2} + |u|_{\alpha}^{2} \right)^{1/2},
$$
and let
$$
H^{\alpha}(\mathbb{R}) = \overline{C_{0}^{\infty}(\mathbb{R})}^{\|.\|_{\alpha}}.
$$

\noindent
We note a function $u\in L^{2}(\mathbb{R})$ belong to $I_{-\infty}^{\alpha}(\mathbb{R})$ if and only if
\begin{equation}\label{FDEeq03}
|w|^{\alpha}\widehat{u} \in L^{2}(\mathbb{R}).
\end{equation}
Especially
\begin{equation}\label{FDEeq04}
|u|_{I_{-\infty}^{\alpha}} = \||w|^{\alpha}\widehat{u}\|_{L^{2}}.
\end{equation}
Therefore $I_{-\infty}^{\alpha}(\mathbb{R})$ and $H^{\alpha}(\mathbb{R})$ are equivalent with equivalent semi-norm and norm. Analogous to $I_{-\infty}^{\alpha}(\mathbb{R})$ we introduce $I_{\infty}^{\alpha}(\mathbb{R})$. Let the semi-norm
$$
|u|_{I_{\infty}^{\alpha}} = \|_{x}D_{\infty}^{\alpha}u\|_{L^{2}}
$$
and norm
\begin{equation}\label{FDEeq05}
\|u\|_{I_{\infty}^{\alpha}} = \left( \|u\|_{L^{2}}^{2} + |u|_{I_{\infty}^{\alpha}}^{2} \right)^{1/2},
\end{equation}
and let
$$
I_{\infty}^{\alpha}(\mathbb{R}) = \overline{C_{0}^{\infty}(\mathbb{R})}^{\|.\|_{I_{\infty}^{\alpha}}}.
$$
Moreover $I_{-\infty}^{\alpha}(\mathbb{R})$ and $I_{\infty}^{\alpha}(\mathbb{R})$ are equivalent , with equivalent semi-norm and norm \cite{VEJR}.

Now we recall the Sobolev lemma.
\begin{Thm}\label{FDEtm01}
\cite{CT-2} If $\alpha > \frac{1}{2}$, then $H^{\alpha}(\mathbb{R}) \subset C(\mathbb{R})$ and there is a constant $C=C_{\alpha}$ such that
\begin{equation}\label{FDEeq06}
\sup_{x\in \mathbb{R}} |u(x)| \leq C \|u\|_{\alpha}
\end{equation}
\end{Thm}
\begin{Remark}\label{FDEnta01}
If $u\in H^{\alpha}(\mathbb{R})$, then $u\in L^{q}(\mathbb{R})$ for all $q\in [2,\infty]$, since
$$
\int_{\mathbb{R}} |u(x)|^{q}dx \leq \|u\|_{\infty}^{q-2}\|u\|_{L^{2}}^{2}
$$
\end{Remark}

The following lemma is a version of the concentration compactness principle

\begin{Lem}\label{FDElm02}
Let $r>0$ and $q\geq 2$. Let $(u_{n}) \in H^{\alpha}(\mathbb{R})$ be bounded. If
\begin{equation}\label{FDECeq01}
\lim_{n\to \infty} \sup_{y\in \mathbb{R}} \int_{y-r}^{y+r} |u_{n}(t)|^{q}dt \to 0
\end{equation}
then $u_{n} \to 0$ in $L^{p}(\mathbb{R})$ for any $p>2$.
\end{Lem}

\noindent
{\bf Proof.} Let $q < s < \beta$ and $u\in H^{\alpha}(\mathbb{R})$, H\"older inequality implies that, if $I_{y} = [y-r, y+r]$,
\begin{eqnarray*}
\|u\|_{L^{s}(I_{y})} &= &\left( \int_{I_{y}} |u(t)|^{s}dt\right)^{1/s} = \left( \int_{I_{y}} |u(t)|^{s(1-\lambda)}|u(t)|^{s\lambda} dt\right)^{1/s}\\
&  \leq &\left( \int_{I_{y}}|u(t)|^{q}dt\right)^{\frac{1-\lambda}{q}} \left( \int_{I_{y}} |u(t)|^{\beta}dt\right)^{\frac{\lambda}{\beta}} \leq \left( \int_{I_{y}}|u(t)|^{q}dt \right)^{\frac{1-\lambda}{q}} \|u\|_{\infty}^{\lambda} |I_{y}|^{\lambda/\beta},
\end{eqnarray*}
where $0< \lambda < 1$ and $\frac{s(1-\lambda)}{q} + \frac{s\lambda}{\beta} = 1$. Now by Theorem \ref{FDEtm01},
$$
\|u\|_{\infty} \leq C \left( \int_{I_{y}}[ |{_{-\infty}}D_{t}^{\alpha}u(t)|^{2} + u(t)^{2}]dt\right)^{1/2}.
$$
So
$$
\|u\|_{L^{s}(\mathbb{R})}^{s} \leq (2r)^{\frac{s\lambda}{\beta}} \left( \int_{I_{y}}|u(t)|^{q}dt \right)^{\frac{(1-\lambda)s}{q}}C^{\lambda s} \left( \int_{I_{y}}[ |{_{-\infty}}D_{t}^{\alpha}u(t)|^{2} + u(t)^{2}]dt\right)^{\frac{\lambda s}{2}}.
$$
Choosing $\lambda s = 2$, i.e $s = 2 + q(1- \frac{2}{\beta})$, wich gives, as $t>s$ is arbitrary, $2< s < 2+ q$, we obtain
$$
\|u\|_{L^{s}(I_{y})}^{s} \leq C' \left(\int_{I_{y}} |u(t)|^{q}dt \right)^{\frac{(1-\lambda)s}{q}} \left( \int_{I_{y}}[ |{_{-\infty}}D_{t}^{\alpha}u(t)|^{2} + u(t)^{2}]dt\right).
$$
Consequently, 
\begin{eqnarray*}
\int_{\mathbb{R}} |u(t)|^{s}dt & = & \sum_{k\in \mathbb{Z}} \int_{2rk}^{2r(k+1)} |u(t)|^{s}dt\\
& \leq & C' \sum_{k\in \mathbb{Z}} \left\{ \left( \int_{2rk}^{2r(k+1)}|u(t)|^{q}dt\right)^{\frac{(1-\lambda)s}{q}} \left( \int_{2rk}^{2r(k+1)}[ |{_{-\infty}}D^{\alpha}u(t)|^{2} + u(t)^{2}]dt\right)\right\} \\
&\leq & \sup_{y\in \mathbb{R}} \left( \int_{y-r}^{y+r} |u(t)|^{q}dt \right)^{\frac{(1-\lambda)s}{q}}\|u\|_{\alpha}.
\end{eqnarray*}
Applying this inequality to each $u_{n}$, we see that $u_{n} \to 0$ in $L^{s}(\mathbb{R})$ for $2<s<q+2$. As $u_{n}\in L^{r}(\mathbb{R})$ for each $r>2$, it follows by interpolation that $u_{n} \to 0$ in $L^{p}(\mathbb{R})$ for each $p>2$. $\Box$

Now we introduce more notations and some necessary definitions. Let $\mathfrak{B}$ be a real Banach space, $I\in C^{1}(\mathfrak{B},\mathbb{R})$, which means that $I$ is a continuously Fr\'echet-differentiable functional defined on $\mathfrak{B}$. Recall that $I\in C^{1}(\mathfrak{B},\mathbb{R})$ is said to satisfy the (PS) condition if any sequence $\{u_{k}\}_{k\in \mathbb{N}} \in \mathfrak{B}$, for which $\{I(u_{k})\}_{k\in \mathbb{N}}$ is bounded and $I'(u_{k}) \to 0$ as $k\to +\infty$, possesses a convergent subsequence in $\mathfrak{B}$.

Moreover, let $B_{r}$ be the open ball in $\mathfrak{B}$ with the radius $r$ and centered at $0$ and $\partial B_{r}$ denote its boundary. We obtain the existence of solutions of (\ref{Eq02}) by use of the following well-known Mountain Pass Theorems, see \cite{PR}.
\begin{Thm}\label{FDEtm02}
Let $\mathfrak{B}$ be a real Banach space and $I\in C^{1}(\mathfrak{B}, \mathbb{R})$ satisfying (PS) condition. Suppose that $I(0) = 0$ and
\begin{itemize}
\item[i.] There are constants $\rho , \beta >0$ such that $I|_{\partial B_{\rho}} \geq \beta$, and
\item[ii.] There is and $e\in \mathfrak{B} \setminus \overline{B_{\rho}}$ such that $I(e)\leq 0$.
\end{itemize}
Then $I$ possesses a critical value $c\geq \beta$. Moreover $c$ can be characterized as
$$
c = \inf_{\gamma \in \Gamma} \max_{s\in [0,1]} I(\gamma (s))
$$
where
$$
\Gamma = \{\gamma \in C([0,1] , \mathfrak{B}):\;\;\gamma (0) = 0,\;\;\gamma (1) = e\}
$$
\end{Thm}
\section{Ground state}

In this section we consider 
\begin{eqnarray}\label{limeq01}
&_{t}D_{\infty}^{\alpha}{_{-\infty}}D_{t}^{\alpha}u(t) + u(t) = f(t,u(t)),\;\;t \in \mathbb{R},\\
& u\in H^{\alpha}(\mathbb{R}),\nonumber
\end{eqnarray}
where $\frac{1}{2} < \alpha <1$ and $f\in C(\mathbb{R} \times \mathbb{R}, \mathbb{R})$ satisfies $(f_{0}) - (f_{4})$

We prove the existence of weak solution of (\ref{limeq01}) by applying the mountain pass theorem to the functional $I$ defined on $H^{\alpha}(\mathbb{R})$ as
\begin{equation}\label{limeq02}
I(u) = \frac{1}{2}\int_{\mathbb{R}}[|{_{-\infty}D_{t}^{\alpha}}u(t)|^{2} + u^{2}(t)]dt - \int_{\mathbb{R}}F(t,u(t))dt.
\end{equation}
However, the direct application of the mountain pass theorem is not sufficient, since the Palais - Smale sequence might lose compactness in the whole space $\mathbb{R}$. To overcome this difficulty, we use a comparison argument devised by Rabinowitz \cite{PR1}, based on the energy functional
\begin{equation}\label{limeq03}
\overline{I}(u) = \frac{1}{2}\int_{\mathbb{R}} [|{_{-\infty}D_{t}^{\alpha}}u(t)|^{2} + u^{2}(t)]dt - \int_{\mathbb{R}}\overline{F}(u(t))dt.
\end{equation}   
It can be proved that the functional $I$ is of class $C^{1}$ and we have
\begin{equation}\label{limeq04}
I'(u)v = \int_{\mathbb{R}} {_{-\infty}}D_{t}^{\alpha}u(t) {_{-\infty}}D_{t}^{\alpha}v(t)dt + u(t)v(t) dt - \int_{\mathbb{R}}f(t,u(t))v(t)dt.
\end{equation}
We notice that a necessary condition for $u\in H^{\alpha}(\mathbb{R}$ to be a critical point of $I$ is $I'(u)u=0$. This condition define the Nehari manifold associated  to the functional $I$ 
$$
\Lambda = \{u\in H^{\alpha}(\mathbb{R})\backslash \{0\}:\;\; I'(u)u=0\}
$$
and so, all non trivial solution of (\ref{limeq01}) belong to the Nehari manifold.

Next, from $(f_{2})$ and $(f_{3})$) it is standard to prove that,  for any $\epsilon >0$, there exists $C_{\epsilon}$ such that
\begin{equation}\label{limeq05}
|f(t, \xi)| \leq \epsilon |\xi| + C_{\epsilon} |\xi|^{p_{0}},\;\;\forall t \in \mathbb{R}
\end{equation}
and consequently 
\begin{equation}\label{limeq06}
|F(t, \xi)| \leq \frac{\epsilon}{2} |\xi|^{2} + \frac{C_{\epsilon}}{p_{0} + 1} |\xi|^{p_{0}+1},\;\;\forall t \in \mathbb{R}.
\end{equation}
We start our analysis with

\begin{Lem}\label{limlm01}
Assume the hypotheses ($f_{0}$)-($f_{4}$) hold. For any $u\in H^{\alpha}(\mathbb{R})\setminus \{0\}$, there is  a unique $\sigma_{u} = t(u) > 0$ such that $\sigma_{u}u \in \Lambda$ and we have
$$
I(\sigma_{u}u) = \max_{t\geq 0} I(tu)
$$
\end{Lem}

\noindent
{\bf Proof.}Let $u\in H^{\alpha}(\mathbb{R})\setminus \{0\}$ and consider the function
$
\psi  :  \mathbb{R}^{+} \to \mathbb{R}$ defined as 
     \begin{eqnarray*} \psi (\sigma) = I(\sigma u) = \frac{\sigma ^{2}}{2}\|u\|_{\alpha}^{2} - \int_{\mathbb{R}}F(t, \sigma u)ddt.
\end{eqnarray*}
Then, by (\ref{limeq06}) we have
\begin{eqnarray*}
\int_{\mathbb{R}}F(t,u)dt & \leq & \frac{C\epsilon}{2}\|u\|_{\alpha}^{2} + \frac{CC_{\epsilon}}{p_{0} + 1}\|u\|_{\alpha}^{p_{0}+1}
\end{eqnarray*}
This implies that
$
\psi (t) > 0,\;\;\mbox{for}\;\;t\;\;\mbox{small}. 
$
On the other hand, by ($f_{1}$) exists $A>0$ such that $F(t, \xi) \geq A|\xi|^{\theta},\;\;\forall |\xi|>1$. So
\begin{eqnarray}\label{limeq07}
I(\sigma u) & \leq & \frac{\sigma^{2}}{2}\|u\|_{\alpha}^{2} - A\sigma^{\theta}\int_{\mathbb{R}}|u(t)|^{\theta}dt
\end{eqnarray}
and since $\theta > 2$, we see that $\psi (\sigma) < 0$ for $\sigma$ large. By ($f_{0}$), $\psi (0) = 0$, therefore there is $\sigma_{u} = \sigma(u) > 0$ such that
$$
\psi (\sigma_{u}) = \max_{\sigma\geq 0} \psi (\sigma) = \max_{\sigma \geq 0} I(\sigma u) = I(\sigma_{u}u).
$$
We see that $\psi'(\sigma) = 0$ is equivalent to
\begin{eqnarray*}
\|u\|_{\alpha}^{2} = \int_{\mathbb{R}}\frac{f(t, \sigma u)u}{\sigma }dx,
\end{eqnarray*}
from where, using ($f_{4}$) we prove that   there is a unique $\sigma_{u} >0$ such that
$
\sigma_{u}u \in \Lambda.
$
$\Box$

Now we define two critical values as follows
\begin{equation}\label{limeq08}
c^{*} = \inf_{u\in \Lambda} I(u)\quad \mbox{and} \quad
c = \inf_{\gamma \in \Gamma} \sup_{\sigma \in [0,1]} I(\gamma (\sigma)),
\end{equation} 
where $\Gamma$ is given by 
$$
\Gamma = \{ \gamma \in C([0,1], H^{\alpha}(\mathbb{R}) /\;\;\gamma (0) = 0,\;I(\gamma (1)) < 0 \}.
$$
Under our assumptions, certainly $\Gamma$ is not empty and $c >0$. The following lemma is crucial and it uses ($f_{1}$).
\begin{Lem}\label{limlm02}
\begin{equation}\label{limeq09}
c^{*} = \inf_{u\in H^{\alpha}(\mathbb{R})\setminus \{0\}} \sup_{\sigma \geq 0} I(\sigma u) = c.
\end{equation}
\end{Lem}

\noindent
{\bf Proof.}
We notice that $I$ is bounded below on $\Lambda$, since by ($f_{1}$), $I(u) >0$, $\forall \; u\in \Lambda$, so that $c^{*}$ is well defined.
By Lemma \ref{limlm01} for any $u\in H^{\alpha}(\mathbb{R})\setminus \{0\}$ there is a unique $\sigma_{u} = \sigma(u) >0$ such that
$\sigma_{u}u \in \Lambda$, then
$$
c^{*} \leq \inf_{u\in H^{\alpha}(\mathbb{R})\setminus \{0\}} \max_{\sigma \geq 0} I(\sigma u).
$$
On the other hand, for any $u\in \Lambda$, we have
\begin{eqnarray*}
I(u) & = & \max_{\sigma \geq 0} I(\sigma u) \geq \inf_{u\in H^{\alpha}(\mathbb{R})\setminus \{0\}} \max_{\sigma \geq 0} I(\sigma u)
\end{eqnarray*}
so
$$
c^{*}=\inf_{\Lambda} I(u) \geq \inf_{u\in H^{\alpha}(\mathbb{R})\setminus \{0\}} \max_{\sigma \geq 0} I(\sigma u),
$$
therefore the first equality in (\ref{limeq09}) holds.
Next we  prove the other equality, that is $c^{*} = c$. We claim that for every $\gamma\in \Gamma$
there exists $\sigma_{0} \in [0,1]$ such that $\gamma (\sigma_{0}) \in \Lambda$.
To prove the claim, we first see that, by (\ref{limeq05}) and the remark \ref{FDEnta01} we have
\begin{equation}\label{limeq10}
\int_{\mathbb{R}} f(t, u)udx \leq \epsilon C \|u\|_{\alpha}^{2} + C_{\epsilon}C \|u\|_{\alpha}^{p_{0}+1}.
\end{equation}
and then, for $\gamma\in\Gamma$ we have 
\begin{eqnarray*}
I'(\gamma (\sigma))\gamma (t) &=& \| \gamma (\sigma) \|_{\alpha}^{2} - \int_{\mathbb{R}} f(t, \gamma (\sigma))\gamma (\sigma)dt
\\& \geq &  \left( 1 - \epsilon C - C_{\epsilon} C \| \gamma (\sigma)\|_{\alpha}^{p_{0}-1} \right) \|\gamma (\sigma)\|_{\alpha}^{2}.
\end{eqnarray*}
If we take $r = \left( \frac{1 - \epsilon C}{C_{\epsilon} C} \right)^{\frac{1}{p_{0}-1}}$, then we see that 
\begin{eqnarray*}
I'(\gamma (\sigma)) \gamma (\sigma ) > 0  \quad \forall \; \sigma \in [0,1], \;\;\mbox{such that,}\;\; \|\gamma (\sigma )\|_{\alpha} < r.
\end{eqnarray*}
On the  other hand, using $(f_1)$ and since $I(\gamma (1)) <0$, we have
\begin{eqnarray*}
\| \gamma (1)\|_{\alpha}^{2} & < & \int_{\mathbb{R}} 2F(t, \gamma (1))dt <  \int_{\mathbb{R}} \theta F(t, \gamma (1))dt \leq  \int_{\mathbb{R}} f(t, \gamma (1)) \gamma (1)dt,
\end{eqnarray*}
that implies 
$
I'(\gamma (1))\gamma (1) < 0.
$
Thus, by the Intermediate Value Theorem, there exists $\sigma_{0}\in (\sigma_{*}, 1)$ such that
$
I'(\gamma (\sigma_{0}))\gamma (\sigma_{0}) = 0$ and so   $\gamma (\sigma_{0}) \in \Lambda.
$, completing the proof of the claim.
From this result,
 $\max_{\sigma \in [0,1]} I(\gamma (\sigma)) \geq I(\gamma (\sigma_{0})) \geq \inf_{\Lambda} I$ and  then
\begin{equation}\label{limeq11}
c\geq c^{*}
\end{equation}

\noindent 
In order to prove the other inequality we see that from (\ref{limeq07}),  there exists $\sigma_{u}^{*}$ large enough such that $I(\sigma_{u}^{*} u) <0$. Now we define the curve
$\gamma_{u}  :  [0,1] \to H^{\alpha}(\mathbb{R})$ as $\gamma_{u}(\sigma) = \sigma (\sigma_{u}^{*}u)$.
Then
$ \gamma_{u} (0)  =  0,$ $I(\gamma (1)) = I(\sigma_{u}^{*}u) <0$ and $\gamma_u$ is continuous, so that  $\gamma_{u} \in \Gamma$. Now, by definition of $\gamma_{u}$,
$$
\max_{\sigma \geq 0} I(\sigma u) \geq \max_{\xi \in [0,1]}I(\gamma_{u}(\xi)),\;\;\forall\;H^{\alpha}(\mathbb{R})\setminus \{0\}
$$
then $c^{*} \geq c$, completing the proof.
$\Box$
\begin{Remark}\label{limnta01}
Since $c = \inf_{\Lambda}I$ and any critical point of $I$ lies on $\Lambda$, if $c$ is a critical value of $I$ then it is the smallest positive critical value of $I$.
\end{Remark}
\begin{Lem}\label{limlm03}
Suppose $\{u_{n}\} \in H^{\alpha}(\mathbb{R}$ and there exists $b>0$ such that
\begin{equation}\label{limeq12}
I(u_{n}) \leq b\;\;\mbox{and}\;\;I'(u_{n}) \to 0
\end{equation}
Then either
\begin{itemize}
\item[(i)] $u_{n} \to 0$ in $H^{\alpha}(\mathbb{R})$, or
\item[(ii)] there is a sequence $(y_{n})\in \mathbb{R}$, and $r, \beta > 0$ such that
$$
\liminf_{n\to \infty} \int_{y_{n} - r}^{y_{n} + r} |u_{n}(x)|^{2}dx > \beta.
$$
\end{itemize}
\end{Lem}

\noindent
{\bf Proof.}
By (\ref{limeq12}) it is standard to check, for $k$ large enough
\begin{equation}\label{limeq13}
b + \|u_{n}\|_{\alpha} \geq I(u_{n}) - \frac{1}{\theta}I'(u_{n})u_{n} \geq \left( \frac{1}{2} - \frac{1}{\theta} \right)\|u_{n}\|_{\alpha}^{2}
\end{equation}
and then $\{u_{n}\}$ is bounded in $H^{\alpha}(\mathbb{R})$.

Suppose ($ii$) is not satisfied, then for any $r>0$, (\ref{FDECeq01}) holds. Consequently by Lemma \ref{FDElm02}
\begin{equation}\label{limeq14}
\|u_{n}\|_{L^{p_{0} + 1}} \to 0.
\end{equation}
Then, noticing that
\begin{equation}\label{limeq15}
I'(u_{n})u_{n} = \| u_{n} \|_{\alpha}^{2} - \int_{\mathbb{R}}f(t, u_{n})u_{n}dt,
\end{equation}
by (\ref{limeq05}) and continuous embedding we have
$$
\int_{\mathbb{R}} f(t, u_{n})u_{n}dt \leq \epsilon C\| u_{n} \|_{\alpha}^{2} + C_{\epsilon}\| u_{n} \|_{L^{p_{0}+1}}^{p_{0}+1}
$$
where $p_{0} + 1 > \theta$. So
\begin{equation}\label{limeq16}
I'(u_{n})u_{n} \geq (1-\epsilon C) \| u_{n} \|_{\alpha}^{2} - C_{\epsilon} \| u_{n} \|_{L^{p_{0}+1}}^{p_{0}+1}.
\end{equation}
Choosing an appropriate $C$ and using (\ref{limeq12}) and (\ref{limeq14}),  we find that 
$u_{n} \to 0\;\;\mbox{in}\;\;H_{\alpha}(\mathbb{R})$, i.e., ($i$) holds.
$\Box$

We define $\overline{\Lambda}$, $\overline{\Gamma}$ and $\overline{c}$, replacing $f$ by $\overline{f}$. The following theorem gives the existence of a solution for the limit problem.

\begin{Thm}\label{limtm01}
$\overline{I}$ has at least one critical point with critical value $\overline{c}$
\end{Thm}

{\bf Proof.} By the Ekeland variational principle (see \cite{JMMW}), there is a sequence $u_{n}$ such that
\begin{equation}\label{limeq17}
\overline{I}(u_{n}) \to \overline{c} \quad \mbox{and} \quad \overline{I}' (u_{n}) \to 0.
\end{equation} 
By (\ref{limeq17}) and ($f_{1}$), given $\epsilon >0$, for sufficiently large $n$,
\begin{eqnarray*}
\left( \frac{1}{2} - \frac{1}{\theta} \right)\|u_{n}\|_{n}^{2} & \leq &\left( \frac{1}{2} - \frac{1}{\theta} \right)\|u_{n}\|_{n}^{2} + \int_{\mathbb{R}} \left[ \frac{1}{\theta}\overline{f}(u_{n})u_{n} - \overline{F}(u_{n})\right]dt \\
& & =  \overline{I}(u_{n}) - \frac{1}{\theta} \overline{I}'(u_{n})u_{n} \leq \|u_{n}\|_{\alpha} + \overline{c} + \epsilon,
\end{eqnarray*}
so that $(u_{n})$ is a bounded sequence. Since $H^{\alpha}(\mathbb{R})$ is reflexive space, there is a subsequence $(u_{n})\in H^{\alpha}(\mathbb{R})$ converging weakly to $u$ in $H^{\alpha}(\mathbb{R})$ and strongly in $L_{loc}^{p}(\mathbb{R})$ for $p\in [2,\infty]$. Thus, for such a subsequence and any $\varphi \in C_{0}^{\infty}(\mathbb{R})$,
$$
\lim_{n\to \infty} \overline{I}'(u_{n})\varphi = \overline{I}'(u)\varphi = 0.
$$ 
If we show that $u\neq 0$, then $\overline{I}'(u) = 0$, and then $\overline{I} (u) \geq \overline{c}$. On the other hand, using ($f_{1}$) again, we see that, for every $r >0$
\begin{eqnarray}\label{limeq18}
\overline{I}(u_{n}) - \frac{1}{2}\overline{I}'(u_n)u_{n} & = \int_{\mathbb{R}} (\frac{1}{2}\overline{f}(u_{n})u_{n} - \overline{F}(u_{n}))dt\nonumber\\
& \geq \int_{-r}^{r} (\frac{1}{2} \overline{f}(u_{n})u_{n} - \overline{F}(u_{n})dt.
\end{eqnarray} 
Since $u_{n} \to u$ in $L_{loc}^{p}(\mathbb{R})$, for any $p\in [2,\infty]$, up to a subsequence,
$$
u_{n}(t) \to u(t)\;\;\mbox{a.e. on} \;\;(-r,r), 
$$ 
and there are $h\in L^{2}(-r,r)$ and $g\in L^{p_{0} + 1}(-r,r)$, such that
$$
|u_{n}(t)| \leq h(t)\quad \mbox{and}\quad |u_{n}(t)| \leq g(x)\quad \mbox{a.e on} \;\;(-r,r).
$$
Moreover by (\ref{limeq05}) and (\ref{limeq06}) we get,
\begin{eqnarray*}
|\overline{f}(u_{n}(t))u_{n}(t)| & \leq \epsilon h(t)^{2} + C_{\epsilon} g(t)^{p_{0} + 1} \in L^{1}(-r,r),\\
|\overline{F}(u_{n}(t))| &\leq \frac{\epsilon}{2}h(t)^{2} + \frac{C_{\epsilon}}{p_{0} + 1} g(t)^{p_{0} + 1} \in L^{1}(-r,r). 
\end{eqnarray*}
By Lebesgue's dominated convergence theorem
$$
\int_{-r}^{r}\overline{f}(u_{n})u_{n}dt \to \int_{-r}^{r}\overline{f}(u)udt \quad \mbox{and} \quad \int_{-r}^{r}\overline{F}(u_{n})dt \to \int_{-r}^{r}\overline{F}(u)dt.
$$
Therefore, since $r$ is arbitrary  
$$
\overline{c} \geq \int_{\mathbb{R}}(\frac{1}{2} \overline{f}(u)u - \overline{F}(u))dt.
$$
Now, since $\overline{I}'(u)u = 0$ and
$$
\overline{I}(u) = \overline{I}(u) - \frac{1}{2}\overline{I}'(u)u = \int_{\mathbb{R}}\left( \frac{1}{2}\overline{f}(u)u - \overline{F}(u) \right)dt,
$$
it follows that $\overline{I}(u) \leq \overline{c}$. 

In order to complete the proof, we just need to show that $u$ is non-trivial. For this purpose by Lemma \ref{limlm03} it is possible to find a sequence $y_{n} \in \mathbb{R}$, $r>0$ and $\beta >0$ such that
$$
\int_{y_{n} - r}^{y_{n} + r} u_{n}^{2}(t)dt > \beta , \;\;\forall n.
$$
Now we define $\overline{u}_{n}(x) = u_{n}(x+y_{n})$. We note that $\overline{u}_{n}$ is bounded in $H^{\alpha}(\mathbb{R})$, and so, up to a subsequence, weakly converges in $H^{\alpha}(\mathbb{R})$ to some $u\in H^{\alpha}(\mathbb{R})$ and strongly in $L^{p}(-r,r)$. But
$$
\int_{-r}^{r}|u(t)|^{p} dt = \lim_{n\to \infty} \int_{-r}^{r}|\overline{u}_{n}(t)|^{2}dt = \lim_{n\to \infty} \int_{y_{n}-r}^{y_{n} + r} |u_{n}(t)|^{p}dt > \beta, 
$$ 
that implies $u\neq 0$. $\Box$

Now we prove our main result
\begin{Thm}\label{limtm02}
$I$ has at least one critical point with critical value $c < \overline{c}$.
\end{Thm}

\noindent 
{\bf Proof.} By definition of  $c$ in \equ{limeq09}, for every sequence  $\{\epsilon_{n}\}$, there exists a sequence of   $\{u_{n}\}$ in $H^{\alpha}(\mathbb{R})$ such that $\| u_{n} \|_{\alpha} = 1$,
\begin{equation}
c \leq \max_{\sigma \geq 0} I(\sigma u_{n}) \leq c + \epsilon_{n} \quad\mbox{and}\quad
\label{limeq19}
\max_{\sigma \geq 0} I(\sigma u_{n}) \to c.
\end{equation}
As in the proof of Lemma \ref{limlm02}, associated with each $u_{n}$, there is a function $\gamma_{n} \in \Gamma$ such that
\begin{equation}\label{limeq20}
\max_{\xi \in [0,1]} I(\gamma_{n}(\xi)) \leq \max_{\sigma \geq 0} I(\sigma u_{n}) \leq c + \epsilon_{n}.
\end{equation}
Now let $X = H^{\alpha}(\mathbb{R})$, $K = [0,1]$, $K_{0} = \{0,1\}$, $M = \Gamma$, $\varphi = \gamma_{n}$ and
$$
c_{1} = \max_{\gamma_{n}(K_{0})}I = 0  < c,
$$
then we can use theorem $4.3$ of \cite{JMMW}, to find a  sequence $\{w_{n}\}$  in $H^{\alpha}(\mathbb{R})$ and $\{\xi_{n}\} \subset [0,1]$ such that $I(w_{n}) \in (c - \epsilon_{n}, c + \epsilon_{n})$,
\begin{equation}\label{limeq21}
\| w_{n} - \gamma_{n}(\xi_{n}) \|_{\alpha} \leq \epsilon_{n}^{1/2}
\quad\mbox{and}\quad
\| I'(w_{n}) \|_{(H^{\alpha})'} \leq \epsilon_{n}^{1/2}.
\end{equation}
Now, since
\begin{equation}\label{limeq22}
I(w_{n}) \to c\;\;\mbox{in}\;\;\mathbb{R}\;\;\mbox{and}\;\;I'(w_{k}) \to 0\;\;\mbox{in}\;\;(H^{\alpha}(\mathbb{R}))',
\end{equation}
as in the proof of the Theorem \ref{limtm01}, we show that $\{w_{n}\}$ is bounded in $H_{\rho}^{\alpha}(\mathbb{R}^{n})$. Moreover up to a subsequence
\begin{equation}\label{limeq23}
w_{n} \rightharpoonup w\;\;\mbox{in}\;\;H^{\alpha}(\mathbb{R})\;\;\mbox{and}\;\;w_{n} \to w\;\;\mbox{in}\;\;L_{loc}^{p}(\mathbb{R}),\;\;2 \leq p \leq \infty,
\end{equation}
where $w$ is weak solution of (\ref{limeq01}). By Lemma \ref{limlm03}, there is a sequence $\{y_{n}\} \subset \mathbb{R}$, $\beta >0$ and $r>0$ such that
\begin{equation}\label{limeq24}
\liminf_{n\to \infty} \int_{y_{n} - r}^{y_{n} + r} w_{n}^{2}dt \geq \beta.
\end{equation}
If $\{y_{n}\}$ contains a bounded subsequence, then (\ref{limeq24}) guarantees that $w \neq 0$ and the results follows. If   $\{y_{k}\}$ is an unbounded sequence. We may assume that, for given $R>0$,
\begin{equation}\label{limeq25}
\lim_{n\to \infty} \int_{-R}^{R} |u_{n}|^{2}dt = 0,
\end{equation}
since the contrary implies that $u\neq 0$ following the same argument as above. In order to complete the proof, we first obtain that
\begin{equation}\label{limeq26}
c < \overline{c}.
\end{equation}
To see this, we use the characterization of $c$ and $\overline{c}$ as in Lemma \ref{limeq02}. Let $\overline{u}$ be a non-trivial critical point of $\overline{I}$ given by Theorem \ref{limtm01} and let
$$
\mathcal{A} = \{x\in \mathbb{R}:\quad f(t,\xi )> \overline{f}(\xi)\;\;\mbox{for all}\;\; \xi >0\}.
$$
Then, by ($f_{5}$) and the fact that $\overline{u}$ is non-zero, there exists $y\in \mathbb{R}$ such that the function $u_{y}$, defined as $u_{y}(t) = u(t + y)$, satisfies 
$$
|\{t\in \mathbb{R}:\;\;|u_{y}(t)|>0\} \cap \mathcal{A}| >0,
$$
where $|.|$ denotes the Lebesgue measure. But then
$$
\overline{c} = \overline{I}(u_{y}) \geq \overline{I}(\sigma u_{y}) > I(\sigma u_{y})\quad\mbox{for all}\;\;\sigma >0.
$$
Choosing $\sigma = \sigma^{*}>0$ such that $I(\sigma^{*}u_{y}) = \sup_{\sigma >0}I(\sigma u_{y})$, we find $\sigma^{*}u_{y} \in \Lambda$ and we conclude that
$$
\overline{c} > I(\sigma^{*}u_{y}) \geq \inf_{u\in \Lambda} I(u) = c,
$$
proving (\ref{limeq26}). Now we see that, for $\sigma \geq 0$, from ($f_{5}$) we have
\begin{eqnarray*}
I(\sigma u_{n}) & = & \overline{I}(\sigma u_{n}) - \int_{\mathbb{R}} ( F(t,u_{n}) - \overline{F}(u_{n}))dt\\
& \geq & \overline{I}(\sigma u_{n}) - \int_{\mathbb{R}} Ca(t)(|\sigma u_{n}(t)|^{2} + |\sigma u_{n}(t)|^{p_{0} + 1})dt.
\end{eqnarray*}
Let $\epsilon >0$, Then, by ($f_{5}$) again, there exists $R>0$ such that
$$
\int_{B(0,R)^{c}} Ca(t)(|\sigma u_{n}(t)|^{2} + |\sigma u_{n}(t)|^{p_{0}+1})dt \leq \epsilon,
$$
for $\sigma$ bounded. Then, by (\ref{limeq25}), 
$$
\lim_{n\to \infty} \int_{B(0, R)} Ca(t)(|\sigma u_{n}(t)|^{2} + |\sigma u_{n}(t)|^{p_{0} + 1})dt = 0.
$$
Choosing $\sigma = \sigma^{*}$ such that $\overline{I}(\sigma^{*}u_{n}) = \max_{\sigma \geq 0} \overline{I}(\sigma u_{n})$, we see that $c \geq \overline{c} - \epsilon$. If $\epsilon>0$ is chosen sufficiently small, this contradicts (\ref{limeq26}). $\Box$


\noindent {\bf Acknowledgements:}
C.T. was  partially supported by MECESUP 0607.

\end{document}